\documentclass[11pt,a4paper,reqno]{amsart}
\usepackage{graphicx}
\usepackage{epsfig}
\usepackage{amsmath}
\usepackage{amsfonts}
\usepackage{amssymb}
\usepackage{amscd}

\begin{document}


\bigskip
\title [ Asymptotic properties ]
{Asymptotic properties of goodness of fit tests based on higher order overlapping spacings }

\author{Sherzod M. Mirakhmedov}

 \maketitle
\begin{center}
{V.I. Romanovskiy Institute of Mathematics, Uzbekistan Academy of Sciences. Tashkent. Uzbekistan.}
 \end{center}
 \begin{center}
{ e-mail: shmirakhmedov@yahoo.com}
\end{center}

\bigskip
\bigskip
\textbf{Abstract}. The paper is devoted to tests for uniformity based on sum-functions of overlapping spacings, where the order of spacings can diverge to infinity as the sample size increases. In particular, it is shown that the asymptotic local power of these tests depends significantly on the asymptotic properties of the counterpart statistics based on disjoint spacings. 

\numberwithin{equation}{section}




\bigskip
\bigskip
{\footnotesize \noindent \textbf{Keywords: exponential distribution, spacings, goodness of fit test, Pitman efficiency, local alternatives, power divergence statistics.} 


\bigskip

\begin{center}
\section{Introduction}. \\
\end{center}

 Assume that $ X_{1},...,X_{n-1} $ are the independent identically distributed random variables (r.v.s),  that form a sample drown from an absolutely continuous distribution whose support is the unit interval [0,1]. The goodness of fit problem consists in testing whether the sample comes from a specified distribution with, say, density \textit{f}. Using a probability integral transformation, this problem can be reduced to testing for uniformity on [0,1]. An important approach to testing this hypothesis is to use spacings-statistics defined below.\\
 \indent To keep future notations simple, we assume that the sample is already in ascending order, i.e. $ X_{1}\leq ...   \leq X_{n-1} .$ Put $ X_{0}=0 $ , $ X_{n}=1 $ and $ X_{k}=1+X_{k-n} $, for $ k>n $ circularly.Two types of spacings are mostly considered in the literature: the overlapping spacings defined as $ D_{k,m}=X_{k+m}-X_{k} \, ,$ $ k=0,1,...,n-1 $ , and the disjoint spacings defined as
 $ D_{k\cdot m,m}=X_{(k+1)\cdot m}-X_{k\cdot m} $ , $ k=0,1,...,N $ ,  where an integer \textit{m} , $ 1\leq m < n ,$  is the order of spacings, $ N=nm^{-1}-1 $ is assumed to be integer, without loss of generality for the asymptotic theory we are dealing with. For brevity, we call these spacings  \textit{m}-spacings. In this paper, we consider the following classes of test statistics: 
 \begin{equation} \label{1,1}
 V_{n,m}=\sum_{k=0}^{n-1}h(nD_{k,m})
\end{equation}
and
\begin{equation} \label{1,2}
 V_{n,m}^{*}=\sum_{k=0}^{N}h(nD_{k\cdot m,m}),
\end{equation}

\noindent where \textit{h} is a real-valued function defined on non-negative half axis. The spacings clearly depend on \textit{n}, 
and the function \textit{h} may also depend on \textit{n} in the cases we consider, but for simplicity of notation we 
suppress the extra subscript \textit{n}. Most common examples of both $ V_{n,m} $ and $ V_{n,m}^{*} $  are the Greenwood statistic where $ h(x)=x^{2} $ , the Moran statistic where $ h(x)=-\log x $ , the entropy type statistic, where  $ h(x)=x\log x $, and the Rao statistic  where  $ h(x)=\vert x-m\vert .$  \\
\indent Given that $V_{n,m}$ and $V_{n,m}^{*}$ are highly dependent r.v.s the large-sample theory is the main avenue for studying such statistics. In order to decide which of the tests based on statistics (1.1) and (1.2) performs better, we study their asymptotic powers and Pitman’s asymptotic relative efficiencies (AREs) in testing of the hypothesis  $ H_{0}:\,  f(x)=1 $, $ 0\leq x \leq 1,$ against the family of contamination alternatives 

\begin{equation} \label{1,3}
H_{1,n}: \\f(x)=1+\delta(n)\textit{l}_{n}(x), \quad  0\leq x \leq 1,
\end{equation}

\noindent where $\delta(n)\rightarrow 0  $  as $ n\rightarrow\infty ,$ and the sequence of functions $ \textit{l}_{n}(x) $   satisfies some smoothness conditions (to be made precise in  Section 2). These alternatives converge to $ H_{0} $  at a rate determined by $\delta(n)  $  , whereas the function  $ \textit{l}_{n}(x) $  defines the "path" along which one deviates from the alternative to the hypothesis.\\
  \indent Statistics (1.1) and (1.2) are symmetric in \textit{m}-spacings. In general, one can consider what have been called "asymmetric" statistics of the form $ h_{0}(nD_{0,m})+...+h_{n-1}(nD_{n-1,m}) $ , or a similar statistic based on the disjoint \textit{m}-spacings, where the real functions $ h_{k} $  are defined on the non-negative half axis. It is known (see, e.g., Kuo and Jammalamadaka,1981) that the asymptotic power of a test based on asymmetric statistics against alternatives at a distance of  $ \delta(n)=n^{-1/2} $ from the null hypothesis will exceed the significance level. However, the asymptotic power of these tests strongly depends on the particular alternative, i.e., the function $l_{n}(x),$ and they perform poorly against any other alternative. Moreover, the asymptotically most powerful test is based on the linear statistic, the asymptotic power of which does not depend on \textit{m}. Therefore the asymmetric tests are not very interesting. We focus on the asymptotic properties of symmetric statistics, paying particular attention to statistics of the form (1.1).\\
  \indent Statistics like (1.1) and (1.2) have been studied in the past. Initially, for testing uniformity the sum-functions of 1-spacings were suggested, see, Pyke (1965), Jammalamadaka and Sethuraman (1975) and references therein. It was shown that the power for such symmetric tests can only be obtained when  $ \delta(n)=n^{-1/4} $. However, as a bonus, one obtains tests whose form does not depend on the alternatives. Later, statistics based on \textit{m}-spacings were shown to be asymptotically more efficient (see, e.g., Del Pino (1979), Cressie (1979), Kuo and Jammalamadaka (1981),Jammalamadaka and Kuo (1984),  Misra and van der Meulen (2001)), with their Pitman efficiency increasing monotonically with \textit{m}, and that among the tests based on symmetric statistics with \textit{m} fixed, the Greenwood statistic  is asymptotically the most efficient one. This led to investigations of what happens when the order of spacings \textit{m}  also increases as a function of the sample size \textit{n}. Tests based on the $ V_{n,m} $ and $ V_{n,m}^{*} $  can discriminate alternatives (1.3) with  $ \delta(n)=(nm)^{-1/4} $  which is somewhat closer to the rate $ \delta(n)=n^{-1/2} $  if  $ m=m(n)\rightarrow\infty $ (see, e.g., Hall (1986), Jammalamadaka et al (1989), Mirakhmedov and Naeem (2008)). Hall (1986) obtained the limiting distribution of $ \bar{V}_{n,m}=\sum_{k=0}^{n-1}h((N+1)D_{k,m}) $   under the alternatives (1.3), where  $ \delta(n)=(nm)^{-1/4} $ , $\textit{l}_{n}(x)=\textit{l}(x)$ , $ \textit{l}(x) $ has five continuous derivatives on [0, 1], and, for some $ r\geq 2 $,   the derivative $ h^{(r+1)}(x) $   exists and is continuous in a neighbourhood of $ x=1 $  , and $ n/m^{r}\rightarrow 0 $ . By utilizing the asymptotic distributions of the statistics $ \bar{V}_{n,m} $ and $ \bar{V}_{n,m}^{*}=\sum_{k=0}^{N}h((N+1)D_{k\cdot m,m}) $ under the alternatives (1.3), where $ \delta(n)=(nm)^{-1/4} $, and the sequence of functions $ \textit{l}_{n}(x) $ \,  satisfies some smoothness  conditions (similar to those given in Section 2), Jammalamadaka et al (1989) studied Pitman ARE of the tests assuming that $ h'''(x) $ is bounded in a neighbourhood of $ x=1 $, and $ m=cn^{p} $ , for some $ c>0 $, $ 0<p<1 $. In particular, they conclude that if we are dealing with spacings of the same order then statistics based on overlapping spacings are considered preferable to the statistics based on disjoint spacings. The asymptotic theory here suggests that the larger \textit{m} values are always better.  \\
 \indent Finally, we note that the asymptotic theory of statistics based on disjoint spacings is much more developed than the one based on overlapping spacings. It includes the CLT, the Edgworth type expansion, the large deviations results, first- and second-order efficiencies, efficiencies in Pitman sense, and also under the Bahadur and the intermediate (somewhat between Pitman and Bahadur settings) approaches, see Zhou and Jammalamadaka (1989), Mirakhmedov and Jammalamadaka (2013), Mirakhmedov et al (2011) and  Mirakhmedov (2010), for details. \\
 \indent In this paper, we focus on the statistic $  V_{n,m} $ , which is a sum of highly dependent r.v.s, since in addition to the spacings themselves being dependent random quantities, the terms form a sequence of \textit{m}-dependent r.v.s. We assume a Lyapunov-type condition that is minimal for the CLT to hold for a sum of \textit{m}-dependent r.v.s.  Additionally, we assume that the function \textit{h} has a continuous on $ (0,\infty) $  derivative. The proof of the main Theorem 2.1 differs from the known proofs of similar statements (see, e.g., Creesie (1976), Kuo and Jammalamadaka (1981), Hall (1986), Jammalamadaka et al (1989), Misra and van der Meulen (2001)). In particular, our results demonstrate the interesting fact (not previously noted) that the statistics based on disjoint spacings play a significant role in determining the asymptotic properties of tests based on overlapping spacings. For instance, the ARE of the statistics based on overlapping spacings is presented as product of three factors, the first of which is the ratio of the AREs of the statistics relative to the corresponding counterparts based on disjoint spacings, the second one is the ratio of the orders of the used spacings, and the  third one is the ratio of the computed under the hypothesis correlation coefficients between  based on disjoint spacings corresponding counterpart statistics and Greenwood statistics (see (2.12)). Next, we apply the general results to a family of tests based on power divergence statistics. In particular, it is shown that, if we consider tests based on overlapping and disjoint spacings of the same order, the statistic based on disjoint spacings requires a 1.5  times larger sample size to be comparable with its counterpart based on overlapping spacings.\\
  \indent  Throughout the paper $ \Phi(x) $ denotes the standard normal distribution function, $ Y_{0},Y_{1},... $ is a sequence of independent r.v.s with the common standard exponential distribution, \, $ Z_{k,m}=Y_{k}+...+Y_{k+m-1} $  , which is a Gamma r.v. with the probability density function
  
 \begin{equation} \label{1,4}
\gamma_{m}(u)=\dfrac{1}{\Gamma(m)}u^{m-1}e^{-u}, \quad u>0 ,
\end{equation}
where $ \Gamma(m) $  is the gamma function. Next, $ P_{i} $ ,$ E_{i} $ and $var_{i}$  stand for the probability, expectation and variance, respectively, computed under the hypothesis $ H_{i}, $ $ i=0,1 $.  All asymptotic statements and limits are  as $ n\rightarrow \infty ,$ $d_{n}\sim b_{n}$ stands for the $d_{n}/b_{n}\rightarrow 1$. 
  
 \begin{center}
\section{Results}. \\
\end{center}
 Consider testing for uniformity versus the sequences of alternative distributions
 
 \begin{equation} \label{2,1}
H_{1}: \\f(x)=1+(nm)^{-1/4}\textit{l}_{n}(x), \quad  0\leq x \leq 1,
\end{equation}
where the function $ \textit{l}_{n}(x) $ is continuous on $ [0,1] $ and there exists a continuously differentiable function $ \textit{l}(x) $
 such that $\sup_{0\leq x\leq 1} \vert \textit{l}_{n}(x)-\textit{l}(x)\vert\rightarrow 0 $ ,\quad $ \int_{0}^{1}\textit{l}(x)dx =0 .$\\
 \indent We will consider  tests based on statistics  (1.1) and (1.2) assuming that the function \textit{h}  is \textit{not linear}. Also, we suppose that we reject $ H_{0}$  for large values of the statistics. We will refer to such a test as an \textit{h}-test. Set  
 
 $$\parallel\textit{l}\parallel_{s}^{s}=\int_{0}^{1}\textit{l}^{s}(x)dx , \quad \tau_{m}=m^{-1}cov(h(Z_{0,m}),Z_{0,m})\, ,\quad \varphi(u)=h(u)-Eh(Z_{0,m})-(u-m)\tau_{m} \, ,$$
 
 \begin{equation} \label{2,2}
 \sigma_{m}^{2}=varh(Z_{0,m}) + 2\sum_{j=1}^{m-1}cov\left( h(Z_{0,m}),h(Z_{j,m})\right)-m^{2}\tau_{m}^{2} \ \ ,\,\ \ \sigma_{m}^{*2}=var\varphi(Z_{0,m})=varh(Z_{0,m})-m\tau_{m}^{2} \, ,
 \end{equation} \, 
 
\begin{equation} \label{2,3} 
    \mu_{m}(h)=corr\left(\varphi(Z_{0,m}),(Z_{0,m}-m)^{2}\right) \, ,  
 \end{equation}
 
\begin{equation} \label{2,4}
\mathbb{A}_{0.n}=Eh(Z_{0,m}) \, , \quad \mathbb{A}_{1.n}=Eh(Z_{0,m})+\frac{\sigma_{m}^{*}\sqrt{m+1}}{\sqrt{2n}}\mu_{m}(h)\parallel\textit{l}\parallel_{2}^{2} \, ,
\end{equation}

We emphasize that by Propositions 2.1 and 2.2 of Mirakhmedov (2025) and Eq. (3.5) and (3.12) of below Section 3 it follows 
\begin{equation} \label{2,5}
E_{i}V_{n,m}=n\mathbb{A}_{i.n}(1+o(1))) \,\, ,\quad var_{i}V_{n,m}=n\sigma_{m}^{2}(1+o(1)))\,,\; i=0,1 .
\end{equation}
Also (see, Mirakhmedov and Jamalamadaka (2013)) 
\begin{equation} \label{2,6}
E_{i}V_{n,m}^{*}=N\mathbb{A}_{i.n}(1+o(1))) \, ,\quad var_{i}V_{n,m}^{*}=N\sigma_{m}^{*2}(1+o(1)))\; , \, i=0,1  .
\end{equation}

\indent \textbf{2.1. Asymptotic power of the tests}.  Set $ g(Z_{0,m})=h(Z_{0,m})-Eh(Z_{0,m})-(Y_{0}-1)m\tau_{m}. $\\

\indent \textbf{Theorem 2.1.} Let $ m=o(n) $  and \\
 (\textit{i}) \, \quad \; \quad \quad $ \dfrac{m^{r-1}}{\sigma_{m}^{r/2}n^{(r-2)/2}} E\vert g(Z_{0,m}) \vert ^{r}\rightarrow 0 $ , some $ r>2 , $\\
 (\textit{ii})  Function  \textit{h}  has a continuous on $ (0,\infty) $  derivative $ h' $  and \\
 
\quad \quad $ \limsup_{n\rightarrow\infty} \dfrac{m^{2}}{\sigma_{m}^{2}\sqrt{n}}\left[E \vert h(Z_{0,m}) h'(Z_{0,m})\vert ^{3}\right] ^{1/3}<\infty \, .$ \\
Then

\[ \sup_{-\infty <x<\infty}\vert P_{i}\lbrace (V_{n,m}-n\mathbb{A}_{i.n})/\sigma_{m}\sqrt{n}<x\rbrace -\Phi(x)\vert \rightarrow 0 \, , \;  i=0,1  \, . \]\\
\\
\indent \textbf{Remark 2.1}. Condition (i) may impose an additional constraint on  \textit{m}. For example, for the Greenwood statistic condition (i) is satisfied if $ m=o(n), $ whereas for  Moran statistic and entropy-type statistic  we assume that $ m=o(n^{(r-2)/2(r-1)}) $ , by choosing  \textit{r} large we can get the rate somewhat closer to $ m=o(n^{1/2}) .$ \\

\indent \textbf{Theorem 2.2.} Under the conditions of Theorem 2.1 one has\\
\noindent 1. A test of size $ \alpha>0 $ based on $ V_{n,m}$ rejects $ H_{0} $  if $ \lbrace V_{n,m}>c_{\alpha}\rbrace $, where $ c_{\alpha}=u_{\alpha}\sigma_{m}\sqrt{n}+nEh(Z_{0,m}) ,$  $ u_{\alpha}=\Phi ^{-1}(1-\alpha) .$ The asymptotic power of this test is $\Phi(e_{m}(h)\parallel\textit{l}\parallel_{2}^{2}-u_{\alpha}) $ where 
\begin{equation} \label{2,7}
e_{m}^{2}(h)=\dfrac{(m+1)\sigma_{m}^{*2}}{2\sigma_{m}^{2}}\mu_{m}^{2}(h) .
\end{equation}\\
2. Within the class of  tests based on $ V_{n,m} $ the asymptotically most powerful (AMP) test of size $ \alpha>0 $  is Greenwood test, i.e., $ h(x)=x^{2} $, whose critical point is $ c_{\alpha} =u_{\alpha}\sqrt{2nm(m+1)(2m+1)/3}+nm(m+1)$   and the asymptotic power
\[ \Phi\left( \sqrt{\frac{3m+1}{2(2m+1)}} \parallel \textit{l}\parallel_{2}^{2}-u_{\alpha}\right)\rightarrow \Phi\left( \dfrac{\sqrt{3}}{2}\parallel \textit{l}\parallel_{2}^{2}-u_{\alpha}\right) \, \, \, if \, \textit{m}\rightarrow\infty \, .\]\\
3. Within the class of tests based on $ V_{n,m} $ the Greenwod test is unique AMP for the fixed \textit{m}, but if $ \textit{m}\rightarrow\infty $  the Greenwod test is no longer the unique AMP. \\

\indent In particular, Theorem 2.2 covers the case when condition (2.5) of Theorem 2 of Hall (1986) does  not hold and  weakens  the conditions on the functions $ \textit{l}(x)$  and $ \textit{h}(x) $. It also extends the results of Section 4 in Kuo and Jammalamadaka (1981) and Section 5 in Misra and van der Meulen (2001) to an increasing \textit{m}.\\
\indent Formula (2.7) for $ e_{m}^{2}(h) $ , which is called the efficacy of the  h-test, can be rewritten as

\[e_{m}^{2}(h)=\dfrac{1}{4m\sigma_{m}^{2}}cov^{2}\left( h(Z_{0,m}),(Z_{0,m}-m-1)^{2}\right),\] 
which coincides with what follows from Theorem 4.2 of Kuo and Jammalamadaka (1981), where the case of fixed \textit{m} is considered. 
 
\indent \textbf{Theorem 2.3.} Assume $m=o(n)$ and for some $ r > 2 $
 $$\dfrac{E\vert \varphi(Z_{0,m}) \vert ^{r}}{\sigma_{m}^{*r/2}N^{(r-2)/2}} \rightarrow 0 . $$
Then

$$ \sup_{-\infty <x<\infty}\vert P_{i}\lbrace (V_{n,m}^{*}-N\mathbb{A}_{i.n}^{*})/\sigma_{m}^{*}\sqrt{N}<x\rbrace -\Phi(x)\vert \rightarrow 0 \, , \quad i=0,1 .$$
\\
A statement similar to Theorem 2.3  has been used (mostly indirectly) in the literature (e.g., Del Pino (1979), Kio and Jammalamadaka (1981), Jamalamadaka et al (1989), Mirakhmedov and Jammalamadaka (2013))  with the restriction that the derivatives $ h' $ and $ h'' $ are continuous on $(0,\infty)$ . In what follows we will use the following statement, which is a consequence of Theorem 2.3 and (2.6)(compare with Theorem 3.1 (i) of Mirakhmedov and Jamalamadaka,2013).\\
  
 \textbf{Corollary 2.1.} Under the conditions of Theorem 2.3 the asymptotic power of the test of size $\alpha >0 $ based on statistic $V_{n,m}^{*}$ is equal to \quad $\Phi(e_{m}^{*}(h)\parallel\textit{l}\parallel_{2}^{2}-u_{\alpha}) $, where   for the efficacy   $ e_{m}^{*2}(h) $  it holds 
\begin{equation} \label{2,9}
e_{m}^{*2}(h)=\dfrac{m+1}{2m}\mu_{m}^{2}(h).
\end{equation}

From (2.7) and (2.8) we see that the function $ \mu_{m}(h) $  plays a significant role in determining the asymptotic nature of the tests based on both overlapping and disjoint spacings. To see its meaning we note that (see, e.g., Del Pino (1979))
\begin{equation} \label{2,9}
 corr_{0}(V_{n,m}^{*}, G_{n,m}^{*2} )=\mu_{m}(h)(1+o(1)), 
\end{equation}
where  $ G_{n,m}^{*2} $  stands for the Greenwood statistics based on disjoint \textit{m}-spacings, i.e., version of (1.2) where $ h(x)=x^{2}.$ 
We note that $ \mu_{m}^{2}(h)\leq1 $  , and $ \mu_{m}^{2}(h)=1 $ for any \textit{m}, only if $ h(x)=x^{2}. $ 
From (2.7) and (2.9) we observe the following interesting phenomenon. The asymptotic power of the tests using $ V_{n,m} $ (based on the overlapping spacings) depends on the correlation between two statistics based on disjoint spacings, specifically, its counterpart $ V_{n,m}^{*} $  and Greenwood statistic  $ G_{n,m}^{*} $. Moreover, in fact, (2.7) and (2.8) together with below Lemma 2.1 allowed us to discover  new properties of Pitman AREs of the tests based on statistics $ V_{n,m} $ and $ V_{n,m}^{*}. $ \\
\indent \textbf{2.2. Pitman efficiency.} Let $ W_{n,m_{1}}^{(1)}(h) $ and $ W_{n,m_{2}}^{(2)}(\psi) $ be two sequences of statistics, each of the form (1.1) or (1.2) with, respectively, tuning functions $h$ and $\psi$, and based on spacings of orders $m_{1}$ and $m_{2}$. Consider tests for uniformity, using these statistics, against the alternatives (2.1). Assume a test of size $ \alpha\in(0.1) $  using $ W_{n,m_{1}}^{(1)}(h) $  based on a sample size \textit{n} and has the power tending to $\beta\in(\alpha,1) $ as $ n\rightarrow\infty $. Let $ k(n)$ be a sequence such that the power of the size $ \alpha  $ test using $ W_{k(n),m_{2}}^{(2)}(\psi) $  also tends to $\beta $ as $ k(n)\rightarrow\infty $. Assume that $ \lim (k(n)/n)$  exist and does not depend on particular choice of $ k(n) $. Then this limit is the Pitman ARE of  $ W_{n,m_{1}}^{(1)}(h) $ wrt $ W_{n,m_{2}}^{(2)}(\psi) $ ( see e.g.,Frazer (1957), p.108), viz.,
 
 $$ PE(W_{n,m_{1}}^{(1)}(h),W_{n,m_{2}}^{(2)}(\psi))=\lim_{n\rightarrow\infty} \dfrac{k(n)}{n}.$$  \\
\indent Let's denote the efficacies of  tests using $ W_{n,m_{1}}^{(1)}(h) $ and $ W_{n,m_{2}}^{(2)}(\psi) $  as $ e_{m_{1}}^{2}(h) $  and $ e_{m_{2}}^{2}(\psi) $, respectively. \\
 \indent \textbf{Lemma 2.1}.Assume that for $ W_{n,m_{1}}^{(1)}(h) $ and $ W_{k(n),m_{2}}^{(2)}(\psi) $ the conditions of Theorem 2.1 or/and Theorem 2.3 are satisfied. Then
     $$PE(  W_{n,m_{1}}^{(1)}(h), W_{n,m_{2}}^{(2)}(\psi)  )=\lim\dfrac{m_{1}e_{m_{1}}^{2}(h)}{m_{2}e_{m_{2}}^{2}(\psi)}. $$

\indent \textbf{Remark 2.2}. We note that statements similar to our Lemma 2.1, which can be encountered in the existing literature, are not accompanied by proofs. The corresponding formulas are, in fact, only applicable to tests using statistics based on spacings of the same order. There is also some confusion regarding the exponent 2 in these formulas. Lemma 2.1 demonstrates that Del Pino (1979, p.1062) was incorrect in asserting that "Rao and Sethuraman use the wrong exponent 2 instead of 4..." \\
\indent \textbf{Corollary 2.2.} Let $ m_{1}(n)$ and $m_{2}(n)$ be sequences such that $ m_{1}(n), m_{2}(n)\rightarrow\infty $  and $ m_{2}(n)=o(m_{1}(n)) $ as $n\rightarrow\infty $. Then any h-test based on either overlapping or disjoint $m_{1}$- spacings is superior to $\psi $-tests based on either overlapping or disjoint $m_{2}$- spacings. \\

In what follows, we will use notation $V_{n,m}(h)$ and $V_{n,m}^{*}(h)$ instead of $V_{n,m}$ and $V_{n,m}^{*}$, $ \sigma_{m}^{2}(h) $ and $ \sigma_{m}^{*2}(h)$ instead of $ \sigma_{m}^{2} $ and $ \sigma_{m}^{*2}$. \\
\indent For ARE of $ V_{n,m}(h) $ wrt $ V_{n,m}^{*}(h) ,$ i.e., when both statistics use the same tuning function and the same order of spacings, by Lemma 2.1, (2.7) and (2.8) we obtain
\begin{equation} \label{2,10}
PE(V_{n,m}(h),V_{n,m}^{*}(h))=\lim \dfrac{m\sigma_{m}^{*2}(h)}{\sigma_{m}^{2}(h)} \geq 1.
\end{equation}
The inequality in (2.10) follows from Theorem 2.3 of Cressie (1979). 

\indent \textbf{Remark 2.3}\footnote{The article appeared in   "Statistics and Probability Letters" 225 (2025) 110461, but it omitted the  sentence “ seen in Mirakhmedov (2025).However, … .” .}. Cressie (1979) primarily deals with fixed \textit{m} and alternatives (1.3), where $ \delta(n)=n^{-1/4} $ and $ \textit{l}_{n}(x)=\phi(x) $, where  $ \phi(x) $   is a bounded positive step function with finite number of discontinuities. Under these conditions, in particular, he proved optimality of Greenwood test within the class of symmetric tests. Note that his formula (2.3) for $B_{m}  $, taken from Holst (1979), needs adjustment, as seen in Mirakhmedov (2024). However, the proof of inequality $ \sigma_{m}^{2}\leq m\sigma_{m}^{*2} $ relies on the fact that  $ \sigma_{m}^{2} $ and  $ \sigma_{m}^{*2} $ are asymptotic values of the variances of $ V_{n,m}(h) $  and $ V_{n,m}^{*}(h) $ , respectively, rather than their exact formula. Therefore, we can  utilize his Theorem 2.3.\\
\indent Equality of (2.10) yields 
 \[PE(V_{n,m}(h),V_{n,m}^{*}(h))=\lim \dfrac{var_{0}(mV_{n,m}^{*}(h))}{var_{0}V_{n,m}(h)} =\lim  \dfrac{var_{0}(N^{-1}V_{n,m}^{*}(h))}{var_{0}(n^{-1}V_{n,m}(h))}  . \]
Note that both  $ N^{-1}V_{n,m}^{*}(h) $ and  $ n^{-1}V_{n,m}(h) $ have the same asymptotic expectation $ Eh( Z_{0,m})$. So, for the Pitman ARE of $ V_{n,m}(h) $ wrt $ V_{n,m}^{*}(h) $  only their asymptotic variances play a key role. \\
\indent It is interesting to note that the ARE of $V_{n,m_{1}}(h)$ wrt $V_{n,m_{2}}^{*}(\psi)$  depends on the function $\psi$ only through  $\mu_{m_{2}}^{2}(\psi)$.  Specifically, by Lemma 2.1 and equations  (2.7), (2.8) and (2.10) we obtain
\begin{equation} \label{2,11}
PE(V_{n,m_{1}}(h),V_{n,m_{2}}^{*}(\psi))=PE(V_{n,m_{1}}(h),V_{n,m_{1}}^{*}(h))\lim \dfrac{(m_{1}+1)\mu_{m_{1}}^{2}(h)}{(m_{2}+1)\mu_{m_{2}}^{2}(\psi)}.
\end{equation}
Further, by Lemma 2.1,equations (2.7) and (2.10)  we obtain 

$$PE(V_{n,m_{1}}(h),V_{n,m_{2}}(\psi))=\lim\left( \dfrac{m_{1}\sigma_{m_{1}}^{*2}(h)}{\sigma_{m_{1}}^{2}(h)}\left( \dfrac{m_{2}\sigma_{m_{2}}^{*2}(\psi)}{\sigma_{m_{2}}^{2}(\psi)}\right)  ^{-1}\dfrac{m_{1}+1}{m_{2}+1}\dfrac{\mu_{m}^{2}(h)}{\mu_{m}^{2}(\psi)}\right) $$
\begin{equation} \label{2,12}
=\dfrac{PE(V_{n,m_{1}}(h),V_{n,m_{1}}^{*}(h))}{PE(V_{n,m_{2}}(\psi),V_{n,m_{2}}^{*}(\psi))}\lim \dfrac{m_{1}+1}{m_{2}+1}\lim\dfrac{\mu_{m}^{2}(h)}{\mu_{m}^{2}(\psi)}.
\end{equation}
Yet, Lemma 2.1 and (2.8) yield
 \begin{equation} \label{2,13} 
 PE(V_{n,m_{1}}^{*}(h),V_{n,m_{2}}^{*}(\psi))=\lim \dfrac{(m_{1}+1)\mu_{m_{1}}^{2}(h)}{(m_{2}+1)\mu_{m_{2}}^{2}(\psi)}.
 \end{equation}
 We would like to emphasize that formulas (2.10)-(2.13) are valid under the conditions of Theorem 2.1 and(or) Theorem 2.3, where the order of spacings can be finite or tend to infinity.
 
\indent \textbf{2.3. Power divergence statistics.}  Set $ \Psi_{d} $, $d\geq-1 ,$  for the family of functions $\psi_{d}(x):=d^{-1}(d+1)^{-1}(x^{d+1}-1) $, $d>-1,$ $ d\neq0$, $ \psi_{0}(x):=x\log x $ and  $ \psi_{-1}(x):=-\log x $,
 where the cases $d=-1$ and $d=0 $ are given by continuity, i.e., by noting that \quad $\lim_{d\rightarrow 0}(x^{d}-1)=\log x $. The power divergence statistics (PDS)(see Cressie and Read (1984)) is defined by 
   
\begin{equation} \label{2,14} 
 \bar{V}_{n,m}(\psi)=\sum_{k=0}^{n-1}\psi(ND_{k,m}), \; \;  \; \bar{V}_{n,m}^{*}(\psi)=\sum_{k=0}^{N-1}\psi(ND_{k\cdot m,m}), 
\end{equation} 
 where $N=n/m$, $\psi(x)\in \Psi_{d}$. We emphasize that due to equalities (34),(35) and (36) of Mirakhmedov and Jammalamadaka (2013)  
 \begin{equation} \label{2,15} 
 \mu_{m}^{2}(\psi)\rightarrow1 ,\quad  as \quad m\rightarrow\infty ,\quad if \;\psi\in\Psi_{d}.
\end{equation}
Note that ARE of tests based on statistics (2.14) are coincide with those based on statistics $V_{n,m}(\psi)$ and $V_{n,m}^{*}(\psi)$, where $\psi \in \Psi_{d}$, because the efficacy remains unchanged under linear transformations of the statistic.\\
\indent Consider the Greenwood statistics $ G_{n,m}^{2} $  and $ G_{n,m}^{*2} $, i.e., the case $d=1$. We have $\mu_{m}(\psi_{1})=1$,  $ \sigma_{m}^{2}(\psi_{1})=2m(m+1)(2m+1)/3 $ and $\sigma_{m}^{*2}(\psi_{1})=2m(m+1) $. Hence,(2.7) and (2.10), respectively, implies 
\begin{equation} \label{2,16}
 e_{m}^{2}(\psi_{1}) \rightarrow\dfrac{3}{4} \quad \text{and} \quad PE( G_{n,m}^{2},G_{n,m}^{*2})=\dfrac{3}{2} \quad \text{if} \quad m\rightarrow\infty .   
\end{equation}
 Let us consider the Moran statistics, viz.,

\[ M_{n,m}=-\sum_{k=0}^{n-1}\log(nD_{k,m})\, \, , \,\quad \, M_{n,m}^{*}=-\sum_{k=0}^{N}\log(nD_{k\cdot m})\, .  \]
For these statistics by Kuo and Jammalamadaka (1981,p.351) and Mirakhmedov and Jamalamadak (2013, pp. 13, 19) we have 
 
\[ \sigma_{m}^{2}(\psi_{-1})= (2m^{2}-2m+1)\zeta(2,m)-2m+1 ,\quad  \sigma_{m}^{*2}(\psi_{-1})=\zeta(2,m)-m^{-1} ,\]
where $\zeta(2,m)$ is the Hurwitz zeta function. We have (see Abramowitz and Stegun 1972, p. 261) $\zeta(2,m)=m^{-1}+(2m^{2})^{-1}-(6m^{3})^{-1}+o(m^{-3}) .$ On using this by (2.7) and (2.10) we obtain

\begin{equation} \label{2,17}
e_{m}^{2}(\psi_{-1})\rightarrow\dfrac{3}{4} \quad \text{and} \quad PE(M_{n,m},M_{n,m}^{*})=\dfrac{3}{2} \,\quad \, \text{as} \, \quad\, m\rightarrow\infty .
\end{equation}

Let's consider the entropy type statistics:

\[ H_{n,m}=\sum_{k=0}^{n-1}(nD_{k,m})\log(nD_{k,m})\, \, , \quad \, \, H_{n,m}^{*}=\sum_{k=0}^{N}(nD_{k,m})\log(nD_{k\cdot m})\, .  \]
Using Eq.(4.8) of Misra and van der Meulen (2001) and Mirakhmedov and Jammalamadaka (2013, pp.13 and 19) in this case we get
\[ \sigma_{m}^{2}(\psi_{0})=(m(m+1)/\sqrt{2})^{2}\zeta(2,m+2)-m(m+1)(2m-1))/4  ,\quad \sigma_{m}^{*2}(\psi_{0})= m(m+1)\zeta(2,m)-m .\]
So by (2.7) and (2.10) after some algebra we obtain 
\begin{equation} \label{2,18}
 e_{m}^{2}(\psi_{0})\rightarrow\dfrac{3}{4} \quad \text{and} \quad  PE(H_{n,m},H_{n,m}^{*})=\dfrac{3}{2} \, \,, \text{as} \, \, m\rightarrow\infty \,.
\end{equation}
Note that in (2.16), (2.17), and (2.18) there are no restrictions on \textit{m}, except that $m=o(n)$.\\
\indent Next, we consider the Kimbal's statistics, denoted as  $K_{n,m}(d)$ and $K_{n,m}^{*}(d)$, respectively, which are versions of (2.14), where $\psi (x)=\psi_{d}(x)$, $d>-1$ and $d\neq0$. Obviously, the derivative $\psi'''_{d}(x)$ is bounded in some neighbourhood of $x=1$. Therefore, based on Corollaries 2.1 and 3.1 and Theorem 4.1 of Jammalamadak et al. (1989), we can see that: if $ m=cn^{p}$, some $c>0$ and $0<p<1,$  then under the alternatives (2.1) one has: \\
\indent  $ N^{-1/2}[ K_{n,m}(d)-nE\psi_{d}(Z_{0,m}/m)]$ is asymptotically normally distributed whose mean is $\psi''_{d}(1)\parallel l\parallel_{2}^{2}/2 $ and variance is $(\psi''_{d}(1))^{2}/3 $,\\
and\\
\indent $ mN^{-1/2}[ K_{n,m}^{*}(d)-NE\psi_{d}(Z_{0,m}/m)]$ is asymptotically normally distributed whose mean is $\psi''_{d}(1)\parallel l\parallel_{2}^{2}/2 $ and variance is $(\psi''_{d}(1))^{2}/2 $.\\
These statements imply that the asymptotic powers of the size $\alpha >0$ tests based on the statistics $K_{n,m}(d)$  and $K_{n,m}^{*}(d)$ are equal to $ \Phi(\sqrt{3}\parallel l\parallel_{2}^{2}/2 -u_{\alpha})$ and $ \Phi(\parallel l\parallel_{2}^{2}/\sqrt{2} -u_{\alpha})$,  respectively. Hence, corresponding  efficacies  are equal, respectively, to 3/4 and 1/2 . Applying these facts by Lemma 2.1 we obtain 
\begin{equation} \label{2,19}  
e_{m}^{*}(\psi_{d})\rightarrow\dfrac{3}{4} \quad and \quad   PE(K_{n,m}(d),K_{n,m}^{*}(d)))=\frac{3}{2} .
\end{equation}

We note that due to Corollary 2.1 of Jammalamadaka et al (1989) for the cases $d=1$ and $d=2$ the condition $m=cn^{p}$ can be removed, since $\psi'''_{d}(x)$ is bounded in the $(0,\infty).$

\indent Assume $m_{1}=c_{1}n^{p_{1}}$ and $m_{2}=c_{2}n^{p_{2}}$, some $ c_{1},  c_{2}>0 $, $ p_{1}, p_{2}\in (0,1)$. 
Since (2.15), (2.16), (2.17) and (2.18) from (2.11) we obtain: for any $h,\psi \in \Psi_{d}$ it holds 
\begin{equation} \label{2,20}
 PE(V_{n,m_{1}}(h),V_{n,m_{2}}^{*}(\psi))= \begin{cases}0 & \text{if} \; p_{1}<p_{2} \\
 
   \dfrac{3}{2}\dfrac{c_{1}}{c_{2}}&\text{if}\; p_{1}=p_{2} \\

 \infty &\text{if}\; p_{1}>p_{2} .
\end{cases} 
\end{equation}

Equalities (2.16)-(2.20) show that in the similar statements of Jammalamadaka et al (1989, pp.363,364) the exponent $ (1+p)^{-1}$ should be removed. Nonetheless, equalities (2.16)-(2.20) serve as examples conforming the claim: if we compare  tests based on overlapping and disjoint spacings of the same order, the overlapping spacings yield more powerful tests. Additionally, a statistics based on disjoint spacings requires a 3/2 times larger sample size to achieve the same precision (i.e., the same significance level and the same asymptotically non-degenerate power) as the test based on overlapping spacings. \\
\indent  Next, on using (2.16)-(2.19) by (2.12) for any $h, \psi \in \Psi_{d}$, $d\geq -1 ,$ we obtain $ PE(V_{n,m}(h), V_{n,m}(\psi))= 1 ,$ and
\begin{equation} \label{2,21}
 PE(V_{n,m_{1}}(h), V_{n,m_{2}}(\psi))= \begin{cases} 0 &\text{if}\; p_{1}<p_{2}\\
   
 \dfrac{c_{1}}{c_{2}}&\text{if}\; p_{1}=p_{2} \\
  
   \infty &\text{if}\; p_{1}>p_{2} .
\end{cases}
\end{equation} 
From (2.13) and (2.15) it follows that for $ PE(V_{n,m_{1}}^{*}(h), V_{n,m_{2}}^{*}(\psi)) $ where $h, \psi \in \Psi_{d}$, $d\geq -1 $ the same properties as in (2.21) hold.  In other words, the asymptotic relative properties of tests based solely on overlapping spacings are identical to those of tests based solely on disjoint spacings.

\begin{center}
\section{Proofs of Theorems}. \\
\end{center}

\indent  Let $ U_{1},...,U_{n-1} $  be a sample from a uniform on [0,1] distribution. Again we assume that the sample is already in ascending order. Set $ U_{0}=0 $, $ U_{n}=1 $, $U_{k}=1+U_{k-n}=1  $  for $ k>n ,$  and $ S_{k,m}=U_{k+m}-U_{k} $ , 
 $ k=0,1,...,n-1 $, the uniform overlapping \textit{m}-spacings. On using notation of Section 2 we denote    
\[  Q_{n,m}=\sum_{k=0}^{n-1}h_{k}(nS_{k,m}), \quad A_{n,m}= \dfrac{1}{n}\sum_{k=0}^{n-1}Eh_{k}(Z_{k,m}) \, ,   \]
\[  C_{n,m}= \dfrac{1}{n}\sum_{k=0}^{n-1}\left[var h_{k}(Z_{k,m})+2\sum_{j=1}^{k+m-1}cov(h_{k}(Z_{k,m}),h_{j}(Z_{j,m}))\right] ,\quad B_{n,m}= \dfrac{1}{n}\sum_{k=0}^{n-1}cov(h_{k}(Z_{k,m}),Z_{k,m})) \, , \]
\[ g_{k,n}(Z_{k,m})=h_{k}(Z_{k,m})-Eh_{k}(Z_{k,m})-B_{n,m}(Y_{k}-1).  \]

The following statement is a consequence of Theorem 2.1 of Mirakhmedov (2024).\\
 
\indent\textbf{ Assertion 3.1}. Assume  $ m=o(n) $, and \\
(i)\quad $C_{n-m,m}= C_{n,m}(1+o(1)) ,  \quad B_{n-m,m}=B_{n,m}(1+o(1)),$ \\ 
$(ii) \quad \quad \quad \quad \dfrac{m^{r-1}}{n^{r/2}(C_{n,m}-B_{n,m}^{2})^{r/2}}\sum_{k=0}^{n-1}E\vert g_{k,n}(Z_{k,m})\vert ^{r}\rightarrow 0 .$\\
 
 \noindent Then  $$ \sup _{-\infty <x<\infty}\vert P\left\lbrace Q_{n,m}<x\sqrt{n(C_{n,m}-B_{n,m}^{2})}+nA_{n,m}\right\rbrace - \Phi(x)\vert\rightarrow 0 .$$\\

 \indent \textbf{Proof of Theorem 2.1}. In the case where the hypotheses hold, Theorem 2.1 immediately follows from Assertion 3.1 because if we set $ h_0 = ... = h_{n-1} = h $, then $ A_{n,m} = Eh(Z_{0,m}) $, $ C_{n,m} - B_{n,m}^2 = \sigma_m^2 $, and condition (i) is obviously fulfilled. \\
\indent Let's consider the case when alternatives (2.1) are valid. Similarly to relation (3.5) of Kuo and Jammalamadaka (1981), we have 
\begin{equation} \label{3,1}
  nD_{k,m}=nS_{k,m}(1-\varepsilon_{_{k,n}})+O_{p}((nm)^{-1/2}),
\end{equation}
where $ \varepsilon_{k,n}=\textit{l}_{n}(\varsigma_{k,n})(nm)^{-1/4} , \, \,\varsigma_{k,n}=(k+0.5)/n \, ,$ $ O_{p}$ is uniform in \textit{k}. \\
\indent Since (3.1), we need to prove that the conditions of Assertion 3.1 are satisfied if we set $ h_{k}(x)=h(x(1-\varepsilon_{k,n})) ,$ $ k=0,1,...,n-1 . $ Set
\begin{equation} \label{3,2}
\tilde{\mathbb{A}}_{1,n}=\dfrac{1}{n}\sum_{k=0}^{n-1}Eh(Z_{k,m}(1-\varepsilon_{k,m})), \quad \mathbb{S}_{1,n}^{2}=C_{n,m}^{(1)}-(B_{n,m}^{(1)})^{2} ,
\end{equation}
where
\begin{equation} \label{3,3}
C_{n,m}^{(1)}=\dfrac{2}{n}\sum_{k=0}^{n-1}\sum_{i=1}^{m-1}cov\left(h(Z_{k,m}(1-\varepsilon_{k,m})),h(Z_{k,m}(1-\varepsilon_{k+i,m}))\right) + \dfrac{1}{n}\sum_{k=0}^{n-1}varh(Z_{k,m}(1-\varepsilon_{k,m})) ,
\end{equation}

\[ B_{n,m}^{(1)}=\dfrac{1}{n}\sum_{k=0}^{n-1}cov(h(Z_{k,m}(1-\varepsilon_{k,m})),Z_{k,m}) .  \]

On changing variables we obtain
\[ \tilde{\mathbb{A}}_{1,n}=\dfrac{1}{n}\sum_{k=0}^{n-1}E\left[ h(Z_{k,m})T_{m}(Z_{k,m},\varepsilon_{k,m})\right] ,  \]
where $ T_{m}(u,v)=(1-v)^{-m}\exp\left\lbrace -uv/(1-v)\right\rbrace . $ On using a Taylor expansion of $ T_{m}(u,v) $  on  $ v $ we get: uniformly in  \textit{u}
\[ T_{m}(u,v)=1+T'_{m}(u,0)v+ T''_{m}(u,0)v^{2}/2+(T'''_{m}(u,0)+o(1))v^{3}/6 \, .\]
Therefore,
\[ \tilde{\mathbb{A}}_{1,n}=Eh(Z_{0,m})+E(h(Z_{0,m})T'_{m}(Z_{0,m},0))\dfrac{1}{n}\sum_{k=0}^{n-1}\varepsilon_{k,m}+ \dfrac{1}{2}E(h(Z_{0,m})T''_{m}(Z_{0,m},0))\dfrac{1}{n}\sum_{k=0}^{n-1}\varepsilon_{k,m}^{2}  \]
\begin{equation} \label{3,4}
+ \dfrac{1}{6}E(h(Z_{0,m})T'''_{m}(Z_{0,m},0))\dfrac{1}{n}\sum_{k=0}^{n-1}\varepsilon_{k,m}^{3} \,  .
\end{equation}

We have $ T'_{m}(u,0)=-(u-m) , \quad T''_{m}(u,0)=u^{2}-2u(m+1)+m(m+1)= (u-m)^{2}-2(u-m)-m $ and  $T'''_{m}(u,0)= -(u-m)^{3}+6(u-m)^{2}+3(m-2)(u-m)-4m .$ Remind notation (2.2). One can observe that

\[ \mu_{m}(h)=corr(\varphi(Z_{0,m}),\mathbb{Z}_{2,m}) , \quad \sigma_{m}^{*2}=var \varphi(Z_{0,m})  ,\]
where $ \mathbb{Z}_{2,m}=T_{m}^{''}(Z_{0,m},0) .$  Set $ \mathbb{Z}_{3,m}=T_{m}^{'''}(Z_{0,m},0) , $ \quad $ \nu_{m}(h)=corr (\varphi(Z_{0,m}),\mathbb{Z}_{3,m})  .$ We have $ E\mathbb{Z}_{j,m}=0 ,$\quad $ E\mathbb{Z}_{j,m}(Z_{0,m}-m)=0  , $ \quad $ j=2,3 \, . $ On using these equalities and the fact that $ n^{-1}(\textit{l}_{n}^{d}(\varsigma_{0,n})+...+\textit{l}_{n}^{d}(\varsigma_{n-1,n}))=\Vert\textit{l}\Vert_{d}^{d}+o(1) $, by (3.4) we obtain

\[ \tilde{\mathbb{A}}_{1,n}=Eh(Z_{0,m})+\dfrac{1}{2}cov(\varphi(Z_{0,m}),\mathbb{Z}_{2,m})\dfrac{1}{n}\sum_{k=0}^{n-1}\varepsilon_{k,n}^{2}- \dfrac{1}{6}cov(\varphi(Z_{0,m}),\mathbb{Z}_{3,m})\dfrac{1}{n}\sum_{k=0}^{n-1}\varepsilon_{k,n}^{3}\]

\[ =Eh(Z_{0,m})+\dfrac{1}{2\sqrt{nm}}cov(\varphi(Z_{0,m}),\mathbb{Z}_{2,m}) (\Vert \textit{l}\Vert _{2}^{2}+o(1))-\dfrac{1}{6(nm)^{3/4}}cov(\varphi(Z_{0,m}),\mathbb{Z}_{3,m}) (\Vert \textit{l}\Vert _{3}^{3}+o(1))\]
\[ =Eh(Z_{0,m})+\dfrac{\sigma_{m}^{*}\sqrt{m+1}}{\sqrt{2n}}\mu_{m}(h)(\Vert \textit{l}\Vert _{2}^{2}+o(1))+c\sigma_{m}^{*}(m/n)^{3/4}\nu_{m}(h)(\Vert \textit{l}\Vert _{3}^{3}+o(1)) \]
\begin{equation} \label{3,5}
=Eh(Z_{0,m})+\dfrac{\sigma_{m}^{*}\sqrt{m+1}}{\sqrt{2n}}\mu_{m}(h)\Vert \textit{l}\Vert _{2}^{2}(1+o(1))=\mathbb{A}_{1,n}(1+o(1)).
\end{equation}
since $ m=o(n) ,$ \, $ var\mathbb{Z}_{2,m}=2m(m+1) $ and $ var\mathbb{Z}_{3,m}=cm^{3}\, . $ 
By arguments similar to those in the proof (3.5), on using $ T_{m}(u,\varepsilon)=1-((u-m)+o(1))\varepsilon $  we obtain
\begin{equation} \label{3,6}
B_{n,m}^{(1)}=cov(h(Z_{0,m}),Z_{0,m})+o((nm)^{-1/4})E(h(Z_{0,m})(Z_{0,m}-m)^{2})
\end{equation}
and
\begin{equation} \label{3,7}
\dfrac{1}{n}\sum_{k=0}^{n-1}varh(Z_{k,m}(1-\varepsilon_{k,m}))=varh(Z_{0,m})+o((nm)^{-1/4})E\vert h^{2}(Z_{0,m})(Z_{0,m}-m)\vert
\end{equation}
In order to simplify calculations, we assume  $ Eh(Z_{0,m})=0 $ up to below eq. (3.13). Set $ \mathbf{Y}_{i}=Y_{k}+...+Y_{k+i-1} \, ,$  $ \mathbf{Y}_{m-i}^{(1)}=Y_{k+i}+...+Y_{k+m-1} \, ,$ $ \mathbf{Y}_{i}^{(2)}=Y_{k+m}+...+Y_{k+m+i-1} \, ,$ $ i= 1,...,m-1\, . $ Note that $ \mathbf{Y}_{s} \, , \mathbf{Y}_{s}^{(1)} \, ,\mathbf{Y}_{s}^{(2)}$ are independent r.v.s with pdf $ \gamma_{s}(u) $ (see (1.4). On changing variables and using equality $ T_{m}(u,\varepsilon)=1-((u-m)+o(1))\varepsilon $ we obtain

\[E[h(Z_{k,m}(1-\varepsilon_{k,m}))h(Z_{k,m}(1-\varepsilon_{k+i,m}))]   \]

\[=E[h(\mathbf{Y}_{i}(1-\varepsilon_{k,m})+ \mathbf{Y}_{m-i}^{(1)} (1-\varepsilon_{k,m}))h(\mathbf{Y}_{m-i}^{(1)}(1-\varepsilon_{k+i,m})+ \mathbf{Y}_{i}^{(2)}(1-\varepsilon_{k+i,m}))]   \]

\[=E[h(\mathbf{Y}_{i}+ \mathbf{Y}_{m-i}^{(1)} )h((\mathbf{Y}_{m-i}^{(1)}+ \mathbf{Y}_{i}^{(2)})(1-\varepsilon_{k+i,m})(1-\varepsilon_{k,m})^{-1})T_{m}(\mathbf{Y}_{i},\varepsilon_{k,n})T_{m}(\mathbf{Y}_{m-i}^{(1)},\varepsilon_{k,n})T_{m}(\mathbf{Y}_{i}^{(2)},\varepsilon_{k,n})]   \]
\[ =E[h(\mathbf{Y}_{i}+ \mathbf{Y}_{m-i}^{(1)} )h((\mathbf{Y}_{m-i}^{(1)}+ \mathbf{Y}_{i}^{(2)})(1-\varepsilon_{k+i,m}) (1-\varepsilon_{k,m})^{-1})\]
\[\cdot(1-((\mathbf{Y}_{i}-i)+o(1))\varepsilon_{k,m}) (1-((\mathbf{Y}_{m-i}^{(1)}-(m-i))+o(1))\varepsilon_{k,m})(1-((\mathbf{Y}_{i}^{(2)}-i)+o(1))\varepsilon_{k,m})]    \]
\[ =E[h(\mathbf{Y}_{i}+ \mathbf{Y}_{m-i}^{(1)} )h((\mathbf{Y}_{m-i}^{(1)}+ \mathbf{Y}_{i}^{(2)})(1-\varepsilon_{k+i,m}) (1-\varepsilon_{k,m})^{-1})]    \]
\[ +E[h(\mathbf{Y}_{i}+ \mathbf{Y}_{m-i}^{(1)} )h((\mathbf{Y}_{m-i}^{(1)}+ \mathbf{Y}_{i}^{(2)})(1-\varepsilon_{k+i,m}) (1-\varepsilon_{k,m})^{-1})(( (\mathbf{Y}_{i}-i)+ (\mathbf{Y}_{m-i}^{(1)}-(m-i))+(\mathbf{Y}_{i}^{(2)}-i) )\varepsilon_{k,m} \] 
\[+ ((\mathbf{Y}_{i}-i)(\mathbf{Y}_{m-i}^{(1)}-(m-i))+(\mathbf{Y}_{i}-i)(\mathbf{Y}_{i}^{(2)}-i)+(\mathbf{Y}_{m-i}^{(1)}-(m-i))(\mathbf{Y}_{i}^{(2)}-i))\varepsilon_{k,m}^{2}  \]
\begin{equation} \label{3,8}
+ (\mathbf{Y}_{i}-i)(\mathbf{Y}_{m-i}^{(1)}-(m-i))(\mathbf{Y}_{i}^{(2)}-i)\varepsilon_{k,m}^{3}) ]:=\nabla_{1}(k,i,n)+\nabla_{2}(k,i,n).
\end{equation}
Here, for simplicity of notation, we have omitted the terms containing factor o(1) because they do not affect the final result. Further, since the function \textit{h} has a continuous on $ (0,\infty) $ derivative $ h' $  we have by Taylor expansion $ h(x(1-\varepsilon))=h(x)-(h'(x)+o(1))x\varepsilon . $ On using this and equality $ (1-\varepsilon_{k+i,n})(1-\varepsilon_{k,n})^{-1}=1-(\varepsilon_{k+i,n}-\varepsilon_{k,n})(1-\varepsilon_{k,n})^{-1}=1-(\varepsilon_{k+i,n}-\varepsilon_{k,n})(1+O(\varepsilon_{k,n})) $  we get 
\[  \nabla_{1}(k,i,n)=E[h(\mathbf{Y}_{i}+\mathbf{Y}_{m-i}^{(1)})h(\mathbf{Y}_{m-i}^{(1)}+\mathbf{Y}_{i}^{(2)})] \]

\begin{equation} \label{3,9}
-E[h(\mathbf{Y}_{i}+\mathbf{Y}_{m-i}^{(1)})(h'(\mathbf{Y}_{m-i}^{(1)}+\mathbf{Y}_{i}^{(2)})+o(1))(\mathbf{Y}_{m-i}^{(1)}+\mathbf{Y}_{i}^{(2)})(\varepsilon_{k+i,n}-\varepsilon_{k,n})(1+O(\varepsilon_{k,n}))].
\end{equation}
Let $ \nabla_{1,2}(k,i,n) $ stands for the second term of (3.9).  Note that due to smoothness of function $ \textit{l}_{n}(x) $ (see (2.1)) and definition of $ \varepsilon_{k,n} $  ( see (3.1)) by mean value theorem we have $(\varepsilon_{k+i,n}-\varepsilon_{k,n})=(nm)^{-1/4}\textit{l}'_{n}(\varsigma)(\varsigma_{k+i,n}-\varsigma_{k,n})=O((nm)^{-1/4}in^{-1}) ,$ where \, $ \varsigma_{k,n}\leq\varsigma\leq\varsigma_{k+i,n} .$  
\, Also the r.v.s $ \mathbf{Y}_{i}+\mathbf{Y}_{m-i}^{(1)} $ and $ \mathbf{Y}_{m-i}^{(1)}+\mathbf{Y}_{i}^{(2)} $ have the same pdf $ \gamma_{m}(u)\, . $ Use these facts to get 

\[   \vert \dfrac{1}{n}\sum_{k=0}^{n-1}\sum_{i=1}^{m-1}\nabla_{1,2}(k,i,n)\vert =c(n^{5}m)^{-1/4}\sum_{i=1}^{m-1}iE\vert h(\mathbf{Y}_{i}+\mathbf{Y}_{m-i}^{(1)})(h'(\mathbf{Y}_{m-i}^{(1)}+\mathbf{Y}_{i}^{(2)})(\mathbf{Y}_{m-i}^{(1)}+\mathbf{Y}_{i}^{(2)})\vert   \]
\[ =o(m^{2}n^{-1/2})(E\vert h(Z_{0,m})\vert^{3}E\vert h'(Z_{0,m})\vert^{3})^{1/3} \, ,  \]

by Holder’s inequality. Apply this and (3.9) to get

\begin{equation} \label{3,10}
\dfrac{1}{n}\sum_{k=0}^{n-1}\sum_{i=1}^{m-1}\nabla_{1}(k,i,n)=\sum_{i=1}^{m-1}Eh(Z_{0,m})h(Z_{i,m})+o(m^{2}n^{-1/2})(E\vert h(Z_{0,m})\vert^{3}E\vert h'(Z_{0,m})\vert^{3})^{1/3} \, .
\end{equation}
Similar, but lengthier algebra gives 

\begin{equation} \label{3,11}
\dfrac{1}{n}\sum_{k=0}^{n-1}\sum_{i=1}^{m-1}\nabla_{2}(k,i,n)=o(m^{4/3}n^{-1/2})(E\vert h(Z_{0,m})\vert^{3})^{2/3}+o(m^{2}n^{-1/2})(E\vert h(Z_{0,m})\vert^{3}E\vert h'(Z_{0,m})\vert^{3})^{1/3} \, .
\end{equation}
Equations (3.8),(3.9),(3.10) and (3.11) imply

\[  C_{n,m}^{(1)}= \sum_{i=1}^{m-1}Eh(Z_{0,m})h(Z_{i,m}) +o(m^{2}n^{-1/2})(E\vert h(Z_{0,m})\vert^{3}E\vert h'(Z_{0,m})\vert^{3})^{1/3} +o(m^{4/3}n^{-1/2})(E\vert h(Z_{0,m})\vert^{3})^{2/3} \, .  \]
This together with (3.2), (3.3), (3.6) and (3.7) yields
\begin{equation} \label{3,12}
\mathbb{S}_{1,n}^{2}=\sigma_{m}^{2}(1+o(1)) .
\end{equation}
Guided by a similar experience with $ \tilde{\mathbb{A}}_{1,n} , $ using equality $ T_{m}(u,\varepsilon)=1-((u-m)+o(1))\varepsilon $ and (3.6), one can observe that  
\[\dfrac{1}{n}\sum_{k=0}^{n-1}E\vert h(Z_{0,m}(1-\varepsilon_{k,n}))-Eh(Z_{0,m})-B_{n,m}^{(1)}(Y_{k}-1)\vert ^{r} =E\vert  g(Z_{0,m})\vert ^{r}(1+o(1))   \]
By this, (3.5), (3.12), Assertion 3.1, where we put $ h_{k}(x)=h(x(1-\varepsilon_{k,n})) $ ,and on using well-known inequality ( see, e.g., Petrov (1975, p.114)
$\vert \Phi(ax+b)-\Phi(x)\vert \leq \vert b\vert (2\pi)^{-1/2}+ \vert a-1 \vert (a\sqrt{2\pi \textit{e}})^{-1} \, , a>0 \,, $ we complete the proof of Theorem 2.1. \\
\indent \textbf{Proof of Theorem 2.2.} Let $ \alpha_{n}(h) $  and $ \beta_{n}(h) $ denote the size and power, respectively, of the \textit{h}-test. Assume $ \alpha_{n}(h)\rightarrow\alpha\in (0,1) $ and set $ \omega_{n}(h)=\sqrt{n}(\mathbb{A}_{1,n}-\mathbb{A}_{0,n})/ \mathbb{S}_{1,n}$ . This quantity $ \omega_{n}(h) $ plays a role of asymptotic shift of statistic $ V_{n,m} $ for the alternatives (2.1). By Theorem 2.1 we obtain $ \beta_{n}(h)=\Phi (\omega_{n}(h) -u_{\alpha}) +o(1) \, ,$ where $ u_{\alpha} =\Phi^{-1}(1-\alpha).$ On the other hand due to (2.4), (2.5) and (3.12)  $ \omega_{n}(h)=\textit{e}_{m}(h)+o(1) \, .$ The first assertion follows. \\
\indent Next, in order of finding the AMP \textit{h}-test, i.e. the one with the maximum efficacy $ e _{m}(h)$ , against the sequence of alternatives (2.1), we note that the efficacy is unaffected by linear transformation of the statistic. As such we may assume without loss of generality that the statistic is already “normalized” by $ \sigma_{m} $  , and then formula (2.7) can be written as $ \textit{e}_{m}^{2}(h)=\parallel \textit{l}\parallel ^{4}_{2}(4m)^{-1}cov^{2}(\varphi(Z_{0,m}),\mathbb{Z}_{2})\, .$ On the other hand by Cauchy—Schwartz inequality $ cov^{2}(\varphi(Z_{0,m}),\mathbb{Z}_{2})\leq var\varphi(Z_{0,m})var\mathbb{Z}_{2} ,  $ with equality iff $ h(x)=(x-m)^{2} .$ Hence, for a fixed \textit{m} the unique AMP test is generated by the Greenwood statistic. The fact that it is no longer unique AMP if $ m=m(n)\rightarrow\infty $   as  $ n\rightarrow\infty $ follows from (2.16)- (2.20), according which PDSs also generate AMP tests. 
  The proof of Theorem 2.2 is completed .\\
\indent \textbf{Proof of Theorem 2.3} is carried out using Corollary 2.1 of Mirakhmedov (2005) and arguments similar to those in the proof (3.5).\\
\indent \textbf{Proof of Corollary 2.1.} Straightforward. 

\indent \textbf{Proof of Lemma 2.1.} Note that $ W_{n,m_{1}}^{(1)}(h) $   and $ W_{k(n),m_{2}}^{(2)}(\psi) $   are used testing for uniformity against the alternatives (2.1), where the "distance" is defined as $ l_{1,n}(x)(nm_{1})^{-1/4} $  and $ l_{2,n}(x)(k(n)m_{2})^{-1/4} $ , respectively, where the functions $ l_{1,n}(x) $ and $ l_{2,n}(x) $ satisfy the smoothness conditions of the alternatives (2.1), provided that $ l_{1,n}(x)(nm_{1})^{-1/4}=l_{2,n}(x)(k(n)m_{2})^{-1/4} . $ 
Hence, $ k(n)/n = (m_{1}/ m_{2})(\parallel l_{2} \parallel_{1}^{4}/ \parallel l_{1}\parallel_{2}^{4})(1+o(1)) .$ On the other hand, according to Theorem 2.2 and Corollary 2.1 the tests will have the same asymptotic power if  
$ e_{m_{1}}^{2}(h)\parallel l_{l} \parallel_{2}^{4} = e_{m_{2}}^{2}(\psi )\parallel l_{2}\parallel_{2}^{4}. $ That is 
$ \parallel l_{2} \parallel_{2}^{4}/ \parallel l_{1}\parallel_{2}^{4}= e_{m_{1}}^{2}(h)/e_{m_{2}}^{2}(\psi) . $ Lemma 2.1 follows.\\

\textbf{Proof of Corollary 2.2} follows by Lemma 2.1, since efficacies (2.7) and (2.8) should be bounded away from zero and infinity.\\

\textbf{Acknowledgements.} The author thanks the referee  for a careful reading of the manuscript and for constructive comments that considerably improved the presentation of the paper.\\

\textbf{Conflict of interests}. The author has no relevant financial or non-financial interests to disclose.

\end{document}